\documentclass[10pt, conference,draftcls,onecolumn]{ieeeconf}
\usepackage{srcltx,pdfsync}
\usepackage{amsmath}
\usepackage{amsfonts}
\usepackage{amssymb}

\usepackage{mathrsfs}
\usepackage{bbm}

\usepackage{graphicx,psfrag}
\usepackage{url}
\usepackage{color}
\usepackage[colorlinks]{hyperref}
\usepackage{gensymb}
\usepackage{wasysym}

\let\theoremstyle\relax

\usepackage[fancythm]{jphmacros2e}
\usepackage{mathtools}

\usepackage{etoolbox}
\newtoggle{paper}
\newtoggle{report}
\toggletrue{report} 

\makeatletter
\let\NAT@parse\undefined
\makeatother
\usepackage[numbers,sort&compress]{natbib}

\newcommand{\N}{\mathbb{N}}
\newcommand{\K}{\mathcal{K}}

\newcommand{\R}{\mbox{$\mathbb{R}$}}

\newcommand{\rd}{\mathbb{R}^d}

\renewcommand{\adj}{\textnormal{Adj}}

\newcommand{\lap}{\textnormal{Lap}}

\usepackage{mathtools}
\mathtoolsset{showonlyrefs=true}

\DeclareMathOperator{\Exp}{\mathbb{E}}

\newlength{\noteWidth}
\long\def\notes#1{\ifinner
          {\footnotesize #1}
          \else
          \marginpar{\parbox[t]{\noteWidth}{\raggedright\footnotesize #1}}
      \fi\typeout{#1}}

\def\rd#1{{\color{red}{#1}}}

\def\pb#1{\notes{pb's comment: \rd{#1}      }}

\graphicspath{{./figures/}{./figures/Plots/}}

\usepackage{xparse}

\ExplSyntaxOn
\cs_new:Nn \pb_prop_gset_bykeys:Nn
{
	\prop_clear:N \l__pb_temp_prop
	\keys_set:nn { pb/propbykey } { #2 }
	\prop_gset_eq:NN #1 \l__pb_temp_prop
}
\cs_generate_variant:Nn \pb_prop_gset_bykeys:Nn { c }

\keys_define:nn { pb/propbykey }
{
	unknown .code:n = \prop_put:Nxn \l__pb_temp_prop { \l_keys_key_tl } { #1 }
}

\cs_generate_variant:Nn \prop_put:Nnn { Nx }

\NewDocumentCommand{\setupcollaborator}{mm}
{
	\prop_new:c { g_collaborator_#1_prop }
	\pb_prop_gset_bykeys:cn { g_collaborator_#1_prop } { #2 }
}

\NewDocumentCommand{\selectcollaborator}{m}
{
	\prop_map_inline:cn { g_collaborator_#1_prop }
	{
		\tl_set:cn { ##1 } { ##2 }
	}
}
\ExplSyntaxOff

\setupcollaborator{PBmac}
{ 
	NAME=Prabir Mac , laptop or mac mini,
	DiCEbibPATH=/Users/pbarooah/Dropbox/DiCElab_bib,
	figPATH=./figures,
}
\setupcollaborator{MHmac}
{
	NAME=Matt Hale using his mac,
	DiCEbibPATH=../../../BIB/DiCElab_bib,
	figPATH=./figures,
}
\setupcollaborator{KPlinux}
{
	NAME=Kendall using her linux deskip,
	DiCEbibPATH=/Users/ZTT/Dropbox (UFL)/BIB/DiCElab_bib,
	ACbibPATH=/Users/ZTT/Dropbox (UFL)/BIB/austin_bib,
	JBbibPATH=//Users/ZTT/Dropbox (UFL)/BIB/bibs,
	TZbibPATH = /Users/ZTT/Dropbox (UFL)/BIB/TTbib,
	figPATH=/Users/ZTT/Dropbox (UFL)/18_Zeng_SingleZoneID_Journal/figures,
}

\selectcollaborator{PBmac} 

\begin{document}
\bstctlcite{IEEEexample:BSTcontrol} 
\title{Differentially Private Smart Metering: \\ Implementation, Analytics, and Billing}

\author{Matthew Hale, Prabir Barooah, Kendall Parker, Kasra Yazdani\\University of Florida\\Gainesville, Florida USA
	\thanks{PB and KP's work are partially supported by NSF via award 1646229 (CPS-ECCS)}}

\maketitle

\begin{abstract}
Smart power grids offer to revolutionize power distribution by sharing granular
power usage data, though this same data sharing can reveal a great deal
about users, and there are serious privacy concerns for customers. In this paper,
we address these concerns using differential privacy. Differential privacy is
a statistical notion of privacy that adds carefully-calibrated
noise to sensitive data before sharing it, and we apply it to provide strong
privacy guarantees to individual customers. Under this implementation, we
quantify the impact of privacy upon smart grid analytics, which are of interest
to utilities, to quantify the loss in accuracy incurred by privacy. Simultaneously,
we quantify the financial impact of privacy through bounding deviations
in customers' utility bills. Simulation results are provided using actual power
usage data, demonstrating the viability of this approach in practice. 
\end{abstract}

\section{Introduction}
In recent years, the emergence of smart grid technologies has led to a great
deal of research in power systems, e.g.,~\cite{tang13,moslehi10,heinen11}. 
The basic idea underlying smart grids is
that smarter sensing technologies can be applied to individual buildings
to give more granular power usage data over time.
Although these technologies are promising, they implicate
significant privacy concerns among users. The power usage data
gathered in smart grids can be aggregated over time for individual users,
and these datasets can be quite revealing. 

Indeed, both the United States Department of Energy (DOE) and the European
Data Protection Supervisor (EDPS) in the European Union have identified major
privacy concerns in smart grids~\cite{doe10,edps12}. Specifically, the EDPS
has reported that smart grid data can ``provide a detailed breakdown of
energy usage over a long period of time, which can show patterns of use''
and that ``[p]rofiles can thus be developed and then applied back
to individual households and individual members of these households''~\cite[Page 15]{edps12}. In addition, the DOE has noted that these usage
patterns ``could reveal personal details about the lives of consumers, such
as their daily schedules''~\cite[Page 2]{doe10}. 

In this paper, we develop a method for preserving users' privacy while
still allowing smart grids to function normally. This method sends
privatized samples of customers' demand data to the utility company,
providing privacy of data in transit and privacy from the utility itself. 
The goal in this privacy implementation is to enable utility companies to 
perform in-house or third-party data analytics without revealing
sensitive customer data. 

Common existing approaches to privacy in smart grids include
battery-based load hiding methods~\cite{efthymiou10,varodayan11},
which require hardware that not all homes have. 
Homomorphic encryption schemes have also been used~\cite{li10,rottondi13}, though
these methods can be computationally demanding and may support
a limited collection of mathematical operations, which restricts
analytics downstream. Developments in~\cite{sankar13} provide information-theoretic
privacy guarantees, though that work requires power consumption and appliance
use to satisfy certain modeling assumptions, which may not always hold. 
In contrast, this work requires no such assumptions. 

Our privacy implementation is built on the framework of differential privacy,
which is a statistical notion of privacy that originates in the database literature~\cite{dwork13}.
 Differential privacy adds carefully-calibrated noise to sensitive data before sharing it,
and it provides strong, rigorous privacy guarantees in several forms.
First, it is immune to post-processing, in that transforming differentially private
data does not weaken its privacy guarantees~\cite{dwork13}. Second, it is robust
to side information, meaning that learning additional information about
data-producing entities does not weaken differential privacy by much~\cite{kasiv08}. 

Originally, differentially privacy was designed to protect sensitive database entries
each time a database is queried~\cite{dwork13}.
Relative to encryption-based approaches, differential privacy is futureproof,
in the sense that its privacy guarantees do not depend upon certain
calculations being infeasible for an adversary, in contrast to encryption techniques
that require updating encryption keys over time. Moreover, differential privacy
can be significantly less computationally demanding than some encryption approaches because
it requires only generating random numbers. 

Differential privacy has been applied in smart power grids previously~\cite{asghar17},
and the most relevant works in the literature are~\cite{zhao14,acs11,eibl17,sandberg15}. 
Work in~\cite{zhao14,acs11}
requires the presence of an external battery, which homes may not have.
Both~\cite{acs11} and~\cite{eibl17} use the infinite divisibility of the Laplace distribution
to have each agent add gamma-distributed noise, leading to differential privacy
for aggregated information but not for individuals;~\cite{acs11} addresses individual
privacy 
by incorporating homomorphic encryption, though this incurs significant computational
expense as discussed above. 
In contrast, 
this paper will provide differential privacy
guarantees to all users without requiring any encryption. 
The developments of~\cite{sandberg15} derive a tradeoff between individuals' privacy
and accuracy of state estimation. In this paper, we derive tradeoffs as well,
but from the perspective of the impact upon aggregate demand analysis and
user billing. 

Given that privacy threats in smart power grids
stem from aggregation over time, we will use a trajectory-level notion
of differential privacy~\cite{leny14,hale18b}. This different form of privacy provides different
privacy guarantees from the database form, and we will elaborate upon
these differences in Section~\ref{sec:privacy}. 
After implementing privacy, a natural
question is how privacy affects the functioning of the grid. 

Aggregate power consumption
is one very common piece of data of interest to utility companies, and
privacy should still allow for accurate aggregate analyses. 
In this work, we quantify the impact of privacy upon both aggregate
power analysis and billing. 
Although we will implement differential privacy by adding
Laplacian noise, we show that these analyses of private data can
be done as though the noise added were Gaussian, thereby unlocking
tools for analysis from the theory of Gaussian stochastic processes.
We then characterize the worst-case error for both total power
consumption at the network level and billing on the per-customer level.
The contributions of this paper are therefore the privacy implementation
itself, together with statistical bounds on the error it induces
in aggregate analytics and customer billing. 

Work in~\cite{danezis11} develops a general-purpose method for paying
bills in a differentially private manner (i.e., a billing method not dedicated to smart grids), though 
our work takes a different approach: under a mandate of grid customer privacy, we quantify
the impact of that privacy upon billing in the smart grid setting.  
To the best of our knowledge, this is the first work to analyze the impact
of differential privacy upon billing in smart grids.

\begin{figure}[t]
\centering
\includegraphics[scale=0.4]{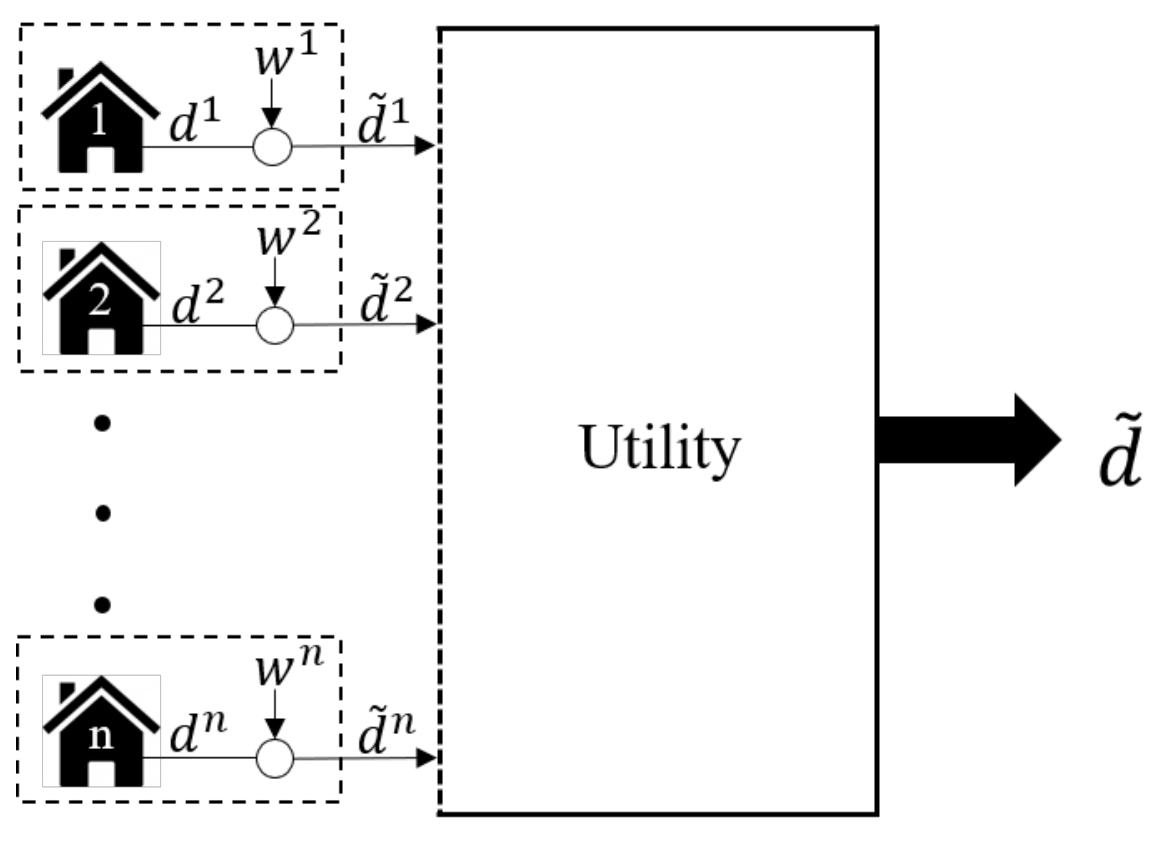}
\caption{In this work, all privacy noise is added at customers' homes before any data
is shared with the utility. Under a mandate of strong user privacy -- even from the utility itself --
we explore the impact of differential privacy upon aggregate load computations and
customer billing.  
}
\label{fig:diagram}
\end{figure}


The rest of the paper is organized as follows. Section~\ref{sec:problem-main}
 states the problems to be solved in this paper. Then,
Section~\ref{sec:privacy} provides the necessary background
on differential privacy and our implementation of it, as well as technical
preliminaries needed in the subsequent analysis. 
Next, Section~\ref{sec:accuracy-analysis} bounds the effects of privacy
upon aggregate analysis and user billing, and Section~\ref{sec:numerics}
provides simulation results to demonstrate these developments on real
smart grid data. 
Finally, Section~\ref{sec:conclusion} provides concluding
remarks and directions for future research.

\section{Problem Description} \label{sec:problem-main}
In this section we outline the problems that are the focus
of the remainder of the paper. 
\subsection{Problem Background and Setup}
While smart meters produce several types of data that can be privatized,
in this paper we focus on privatizing signals of demand data.
Time is measured by a discrete counter, $k$, which represents
sampling continuous-time demand signals
with a sampling period $\Delta t$. In this paper we consider time indices
$k \in \mathcal{K} := \{1, \ldots, K\}$ for some $K \in \N$. 
The total number of customers in the region under study is 
denoted by $N \in \N$. Define $d_k^i$ as the power demand (in kW) 
of the $i^{th}$ customer at time instant $k$, for $i \in [N] := \{1, \ldots, N\}$. 

As noted above, we will protect the power demand data of individual
customers using differential privacy, and that will require noise to be
added to customers' demand signals.
The \emph{private demand} of customer $i$ at time~$k$ is denoted by 
\begin{align}\label{eq:private-demand}
\hat{d}_k^i \eqdef  d_k^i + w_k^i, 
\end{align}
where $w_k^i$ is the noise added to provide differential privacy to customer~$i$. 
The precise distribution of noise
will be developed below in Section~\ref{sec:privacy}, 
and here we only introduce the terminology required to state the problems we solve.

We consider a privacy architecture in which noise is added to the data in the smart meter itself. 
The raw data never leaves the consumer's premises, as shown
in Figure \ref{fig:diagram}. The private data set $\{\hat{d}_k^i\}_{i \in [N], k \in \mathcal{K}}$ is available to the utility, which it can release to an in-house or third-party    
vendor for performing data analytics. 
This architecture provides strong privacy guarantees to customers because not even
the utility company has access to customers' raw power demand data. Motivation for considering this scenario comes from the possibility of mandated privacy rules on smart meter data in the future, much like the European Union's recent General Data Protection Regulation 
(GDPR)\footnote{https://ec.europa.eu/commission/priorities/justice-and-fundamental-rights/data-protection/2018-reform-eu-data-protection-rules\_en} for personal data. 
Motivation also comes from preventing privacy breaches of sensitive customer data,
and, by not sharing raw sensitive data at all, this architecture precludes such breaches. 
For the customer,  a potential cost for such strong privacy  is uncertainty in billing. 
By only sharing privatized demand data, 
bills must be assessed using noisy data. In fact any data analytics performed with private data will 
have some amount of error  
due to privacy noise. We will explore the impact of privacy upon analytics and
billing below.

\subsection{Problem Statements}\label{sec:problem}
We focus on two use cases in particular: (i) energy use of a single
customer and associated billing considerations, and (ii) aggregate demand
of a collection of consumers.

\def\K{\ensuremath{\mathcal{K}}}

\paragraph{Billing/Energy Use} The main variable in the energy bill a customer receives at the end of a billing period is the total energy used in that period\footnote{For industrial customers there are other important variables such as peak demand, but we limit to residential customers here.}. Let the energy consumed by customer $i$ in the $k$-th sampling interval be $E_k^i$. Let the \emph{total energy consumed} by customer $i$ in an interval $\K = \{0,1,\dots,K\}$ be $E^i_\K \eqdef \sum_{k=1}^K E_k^{i}$.  Assuming $\Delta t $ is small enough, $E_k^i = d_k^i \Delta t$, and thus $E^i_\K = \sum_{k \in \K} (d_k^i \Delta t)$. 
However, the utility does not have access to this value exactly. 
Since only private energy data is available to the utility, it instead must bill its customers based on an estimate of the energy computed from the private energy use data. The minimum variance estimator of the total energy consumption of customer $i$ over the interval $\K$, which we call \emph{private energy consumption} during $\K$, is $\hat{E}^i_\K \eqdef \sum_{k=1}^K \hat{d}_k^{i}\Delta t = \sum_{k \in \K} (d_k^i + w_k^i)\Delta t$. The estimation error  is 
\begin{align}\label{eq:error-energy}
\tilde{E}^{i}_\K \eqdef \hat{E}^i_\K - E^i_\K =  \Delta t\sum_{k\in \K} w_k^i .
\end{align}
A statistical characterization of the error $\tilde{E}^{i}_\K$ as a function of privacy parameters is desired in order to understand privacy's impact
upon billing, and that constitutes the first problem we will solve. 

\paragraph{Error in the aggregate demand estimate} The \emph{aggregate demand} $d_k$ of a collection of consumers is the sum of demands of individual households at time~$k$: 
$d_k \eqdef \sum_{i=1}^N d_k^i$. 
The \emph{private aggregate demand}  $\hat{d}_k$ is the utility's 
estimate obtained by using the private data:  $\hat{d}_k \eqdef \sum_{i=1}^N \hat{d}_k^i = \sum_{i=1}^N (d_k^i + w_k^i)$. The error in this estimate is therefore 
\begin{align}\label{eq:error-aggrdemand}
\tilde{d}_k \eqdef \hat{d}_k - d_k =  \sum_{i \in [N]} w_k^i \defeq w_k,
\end{align}
which follows from Equation~\eqref{eq:private-demand}. The second problem 
is to statistically characterize this error and understand its impact upon the utility's analytics
as a function of customers' privacy parameters.

\section{Privacy Implementation} \label{sec:privacy}
In this section we briefly provide the necessary background on
trajectory-level differential privacy, for which a complete exposition
can be found in~\cite{leny14}. 
Then we describe our privacy
implementations for $d^i_k$. 
We emphasize that our developments
are not on the theory of differential privacy itself, but, instead, are
on the application of differential privacy to smart meter data. 

\subsection{Differential Privacy Background}
This subsection provides the requisite differential
privacy background. We retain the notation typically used
in the literature and denote an arbitrary input by~$u$ in this subsection.
These tools will be applied to $d^i_k$  specifically
below. 

We first define the space $\tilde{\ell}_p$ that will contain
users' sensitive input signals. User~$i$'s input is $u^i := \{u^i_k\}_{k \in \mathcal{K}}$. 
The truncation operator $P_T : \tilde{\ell_p} \to \ell_p$ is defined as
\begin{equation}
P_T[u^i] = \begin{cases} u^i_k & k \leq T \\ 0 & k > T \end{cases}.
\end{equation}
Then we define
$\tilde{\ell}_p := \{u^i \mid P_T[u^i] \in \ell_p \textnormal{ for all } T \in \N\}$,
and we have $u^i \in \tilde{\ell}_p$.
In particular, $\ell_p \subsetneq \tilde{\ell}_p$, and the space
$\tilde{\ell}_p$ includes, for example, signals that do not vanish asymptotically.

For all $i \in [N]$, user $i$ contributes the input signal
$u^i \in \tilde{\ell}_p$ to some system. 
The goal of differential privacy is to make ``nearby'' input signals produce
outputs with similar probability distributions. The notion of ``nearby''
is formalized through an adjacency relation, which we define next. 
\begin{definition}
	Fix an adjacency parameter $b > 0$. For two inputs $u^i_1, u^i_2 \in \tilde{\ell}_p$, 
	we define the binary symmetric adjacency relation
	${\adj_b : \tilde{\ell}_p \times \tilde{\ell}_p \to \{0, 1\}}$ via
	\begin{equation}
	\adj_b(u^i_1, u^i_2) = \begin{cases} 1 & \|u^i_1 - u^i_2\|_{\ell_p} \leq b \\ 0 & \textnormal{otherwise} \end{cases}. 
	\end{equation}
\end{definition}
Differential privacy masks the differences between adjacent inputs. For example,
in a smart grid, two demand signals may differ if a home's occupants return from work at
a different time, or if they are on vacation for some period of time before resuming
normal power usage.  
The parameter~$b$ can then be chosen based on the size of difference
that should be masked, and larger values of~$b$ give greater privacy
guarantees because more trajectories must be made approximately
indistinguishable. 
We note that this choice of adjacency pertains to 
single customers, which differs from adjacency relations
used for trajectories in~\cite{leny14,hale18}. 

Differential privacy concerns the probability distributions of outputs that
correspond to adjacent inputs. To state its definition, 
we define a $\sigma$-algebra $\Sigma_p$ over $\tilde{\ell}_p$, and details for
cosntructing $\Sigma_p$ are in~\cite{leny14}.
Differential privacy itself is enforced
by a mechanism. Fixing a probability space $(\Omega, \mathscr{B}, \mathbb{P})$, we formally
state the definition of a differentially private mechanism below. 

\begin{definition} \label{def:dp}
	A mechanism $M : \tilde{\ell}_p \times \Omega \to \tilde{\ell}_q$ preserves
	$\epsilon$-differential privacy if, for all adjacent inputs $u^i_1, u^i_2 \in \tilde{\ell}_p$,
	\begin{equation}
	\mathbb{P}[M(u^i_1) \in S] \leq e^{\epsilon} \mathbb{P}[M(u^i_2) \in S] 
	\textnormal{ for all } S \in \Sigma_p. 
	\end{equation}
\end{definition}

This definition encodes the fact that any eavesdropper or adversary
is unlikely to learn anything meaningful about sensitive data by looking
at differentially private information. 
The likelihood of such privacy breaches is controlled by the privacy
parameter~$\epsilon$, which is user-specified. 
Typical values 
of $\epsilon$ in the literature range from $0.1$ to $\ln 3$~\cite{leny14}, with smaller
values providing stronger privacy guarantees. 
One of the most widely used mechanisms for enforcing Definition~\ref{def:dp}
is the Laplace mechanism, which adds Laplace noise to sensitive data. To define it, we first
define the sensitivity of a system. Regarding a causal, discrete-time,
deterministic dynamical system as a map $\mathcal{G}_i : \tilde{\ell}_p \to
\tilde{\ell}_q$, we have the following definition.

%

\begin{definition} \label{def:sensitivity}
	Let $u^i_1, u^i_2 \in \tilde{\ell}_p$ be adjacent inputs to a system $\mathcal{G}_i$.
	Then the $\ell_p$-sensitivity of $\mathcal{G}_i$ is 
	\begin{equation}
	\Delta_p \mathcal{G}_i := \sup_{\substack{u^i_1, u^i_2 \in \tilde{\ell}_p \\ \adj_b(u^i_1, u^i_2) = 1}}
	\|\mathcal{G}_i(u^i_1) - \mathcal{G}_i(u^i_2)\|_{\ell_p}. 
	\end{equation}
\end{definition}

It is in terms of the $\ell_1$-sensitivity that we will define the Laplace mechanism.
First, we introduce the notation $\lap(c)$ to denote a Laplace distribution
with mean zero and scale parameter $c$, i.e.,
\begin{equation}
\lap(c) := f(0, c; x) = \frac{1}{2c}\exp\left(-\frac{|x|}{c}\right). 
\end{equation}

\begin{definition} \label{def:dp}
	Let $\mathcal{G}_i : \tilde{\ell}_p \to \tilde{\ell}_q$ be a causal, discrete-time, deterministic dynamical system
	with $\ell_1$-sensitivity $\Delta_1 \mathcal{G}_i$. Then the Laplace mechanism
	$M_i(u^i) := \mathcal{G}_i(u^i) + v^i$
	is $\epsilon$-differentially private if $v^i_k \sim \lap(c)$ with
	$c \geq \Delta_1\mathcal{G}_i/\epsilon$. 
\end{definition}



%

We emphasize that privacy here is enforced at the trajectory level and
customer~$i$'s mechanism~$M_i$ provides a single
$\tilde{\ell}_p$-valued query that is shared pointwise in time.
Unlike other works, these trajectory-level privacy guarantees
do not weaken over time because the output at each point in time
is not a query of the initial state. Instead, the private output at 
each point in time merely contributes to assembling a single
trajectory-valued query. 

\subsection{Privatizing Customers' Demand Data}
We consider customers' power demand trajectories over time to be the sensitive
inputs that need to be protected. Thus, while above we have used
the symbol $u$ to represent a sensitive input as in the private control literature, 
we turn now to protecting the trajectories $\{d^i_k\}_{k \in \mathcal{K}}$ for each $i \in [N]$;~\cite[Lemma 2]{leny14} 
guarantees that privacy at the
trajectory level also provides privacy to all finite truncations of trajectories, 
and thus the above definitions apply regardless of the value of~$K$. 

We regard customer~$i$'s demand signal~$d^i$ as passing through
a memoryless ``identity system'' whose output is equal to its input. Formally,
$\mathcal{I}\big(d^i_k\big) = d^i_k$,
for all $k \in \mathcal{K}$ and $i \in [N]$. 
We use the adjacency relation $\adj_b$ for each 
system of this kind, and implementing privacy for it then requires a sensitivity
bound. This bound takes an elementary form in the following lemma.

\begin{lemma} \label{lem:sensitivity}
The identity system $\mathcal{I}(d^i_k) = d^i_k$ has $\ell_p$-sensitivity
bounded by~$b$, i.e., $\Delta_p \mathcal{I} \leq b$.
\end{lemma}
\emph{Proof:} For adjacent~$u^i_1$ and $u^i_2$,  
$\|\mathcal{I}\big(u^i_1\big) - \mathcal{I}\big(u^i_2\big)\|_{\ell_p}
= \|u^i_1 - u^i_2\|_{\ell_p} \leq b$,
which follows from adjacency.
\hfill $\blacksquare$

The basic form of the sensitivity bound above also gives a simple
differential privacy implementation. 
\begin{theorem} \label{thm:lapmech}
Let $\epsilon > 0$ and $b > 0$ be given. The following  Laplace mechanism provides $\epsilon$-differential privacy to the $i$-th customer:
$M_i(d^i) = d^i + w^i, \quad w^i_k \sim \lap\big(0, b/\epsilon)$. 
\end{theorem}
\emph{Proof:} Follows from Definitions~\ref{def:sensitivity}
and~\ref{def:dp} and Lemma~\ref{lem:sensitivity}. \hfill $\blacksquare$

It is understood above that the privacy noise terms added to each customer's data are
mutually independent. 
Since demand data is privatized at each customer's site, the utility has access to only the private data. Any analytics are performed not upon agent~$i$'s
raw demand value $d^i_k$, but instead upon the private demand $\hat{d}_k^i$. 
The privacy noise impacts the accuracy of analytics performed with this data, and quantifying these impacts is the subject
of~Section~\ref{sec:accuracy-analysis}. Before doing so, we derive several
results we will need in our accuracy analyses. 

\subsection{Technical Preliminaries for Stochastic Processes}
As both aggregate demand of a neighborhood and energy use of a single customer involve summing over noisy samples, sums of random variables will appear in the forthcoming analysis.  In either case, as long as the sum is over a large number of random variables, we can appeal to the central limit theorem to use a Gaussian approximation. In case of private aggregate demand analysis, the sum will be over consumers while in case of energy use the sum will be over time. The number of summands in either case is finite, and in some cases may be small. To justify a Gaussian analysis, we now introduce a theoretical tool, the Berry-Esseen
theorem, which provides a convergence rate for the central limit theorem as a function of the number of 
summands. 

\begin{lemma}(\emph{Berry-Esseen Theorem}~\cite[Chapter XVI]{feller08}) \label{lem:berry}
	Let $\{X_{\ell}\}_{\ell \in \mathcal{L} \subseteq \N}$ be a collection
	of~$m$ i.i.d. random variables, with zero mean, variance $\sigma^2$, and third
	moment $\rho < \infty$. Let
	\begin{equation}
	Z(m) = \frac{X_1 + \cdots + X_m}{\sigma \sqrt{m}}
	\end{equation}
	have	CDF $F_{Z(m)}$ and let $\Phi(x)$ be the CDF of a standard normal random
	variable (zero mean, unit variance). Then 
	\begin{equation}
	\sup_{x \in \R} \left|F_{Z(m)}(x) - \Phi(x)\right| \leq C\frac{\rho}{\sigma^3\sqrt{m}},
	\end{equation}
	where $C$ is a universal constant. 
\end{lemma}
Over time, various estimates have been made for~$C$, and recent
work has bounded it above by $0.4748$ for this setting~\cite{shevtsova11}. For simplicity,
we will proceed with $C = \frac{1}{2}$. 

In the problems we consider, 
each~$X_{\ell}$ is a Laplacian random variable with mean 0 and scale parameter  $c$. Specializing to this case, we have the following
result. 
\begin{lemma} \label{lem:minimum-m}
	Let $\eta > 0$ be arbitrary. The distribution of the sum of~$m$ i.i.d. Laplacian random variables, each with distribution~$\lap(c)$, is approximately Gaussian with variance 
	$2m c^2$, with error less than $\eta$ (in the sense of Lemma~\ref{lem:berry})
	if 
		\begin{equation}
	m \geq \frac{9}{64\eta^2 c^6}.
	\end{equation}
\end{lemma}
\emph{Proof: } Straightforward calculations show that the Laplace distribution
$\lap(c)$ has variance $2c^2$ and third moment $6c^3$. 
With $C = \frac{1}{2}$, using
these values in Lemma~\ref{lem:berry} and bounding the error above
by $\eta$ gives $\frac{1}{2}\frac{6c^3}{8c^6\sqrt{m}} \leq \eta$.  Solving for~$m$
completes the proof. \hfill $\blacksquare$

Below we use the term ``$\eta$-approximately Gaussian''
to describe any distribution with error not
more than $\eta$ in the sense of Lemma~\ref{lem:berry}. 
As a concrete example, for $b=2$, $\epsilon = 1$, Lemma \ref{lem:minimum-m} tells us that we need $m \geq 22$ for the sum of $m$ Laplace random variables to be $0.01$-approximately Gaussian. We will later use this bound to show that Gaussian analyses
are possible under reasonable assumptions. 


\begin{lemma} \label{lem:maxgauss}
Let $g_1,\dots,g_n$ be i.i.d. Gaussian, each with 0 mean and variance $\sigma_0^2$.	For $n \geq 43$,
	\begin{equation}
	0.338 \sigma_0 \sqrt{\log n} \leq 
	\mathbb{E}\left[\max_{1 \leq i \leq n} g_i \right]
	\leq \sqrt{2} \sigma_0 \sqrt{\log n}.
	\end{equation}
	The variance of the max, for any $n$, is upper bounded by
	\begin{align}
	\textnormal{var}\left[\max_{1 \leq i \leq n} g_i\right] \leq 4\sigma_0^2.
	\end{align}
\end{lemma}
Although the maximum of Gaussian random variables is widely studied, bounds on the mean available in the literature frequently contain unknown constants, e.g.,~\cite[Appendix A]{ChatterjeeSuperconcentration_book:2014}. The lower bound on the mean provided in the Lemma, which we prove below, represents a slight tightening of the best available bound we have found~\cite{kamath15}. This bound will be applied below to bound the maximum error over samples
in time and the maximum error over customers. Requiring $n \geq 43$ is a very mild assumption
because it merely requires having at least $43$ customers in a power grid or at least
$43$ samples in a dataset, and both of these conditions are easily satisfied in the vast
majority of cases. 
The bound on the variance is derived from standard concentration inequalities on the maximum of Gaussians, and is likely to be already known, but we were unable to find a reference and thus include the proof here. 
\begin{proof-lemma}{\ref{lem:maxgauss}}
The mean upper bound can be found in~\cite[Appendix A]{ChatterjeeSuperconcentration_book:2014}, so we prove only a lower bound for the mean. 
Let $S_1$ be the event ``there exists an~$i$ such
that $g_i \geq \sigma_0\sqrt{\log n}$'', and
let $S_2$ be the complementary event
``$g_i < \sigma_0\sqrt{\log n}$ for all~$i$''. 
For economy of notation we define the symbol $g_{max} = \max_{1 \leq i \leq n} g_i$.
Then
\begin{equation} \label{eq:Cbound1}
\mathbb{E}\left[g_{max}\right] = \mathbb{E}\left[g_{max}\mid S_1\right]\Prob[S_1] + 
\mathbb{E}\left[g_{max} \mid S_2\right]\Prob[S_2].
\end{equation}
We lower-bound the right-hand side by noting that
\begin{align}
\mathbb{E}\left[g_{max} \mid S_2\right] &\geq \mathbb{E}\left[g_i \mid S_2\right] 
\geq \mathbb{E}[g_i \mid g_i < 0] = -\sigma_0\sqrt{2/\pi},
\end{align}
which uses the expectation of the half-normal
distribution~\cite[Equation (3)]{leone61}.

Using the complementarity of $S_1$ and $S_2$, we
return to Equation~\eqref{eq:Cbound1} to find
\begin{equation} \label{eq:Cbound2}
\mathbb{E}\left[g_{max}\right] \geq
\mathbb{E}\left[g_{max} \mid S_1\right]\Prob[S_1] - \sigma_0\sqrt{2/\pi}\big(1 - \Prob[S_1]\big).
\end{equation}
Certainly
$\mathbb{E}\left[g_{max} \mid S_1\right] \geq 
\sigma_0\sqrt{\log n}$,
which gives
\begin{equation} \label{eq:Cbound4}
\mathbb{E}\left[g_{max}\right]
\geq \big(\sigma_0\sqrt{\log n}\big)\Prob[S_1]
-\sigma_0\sqrt{2/\pi}\big(1 - \Prob[S_1]\big),
\end{equation}
and all that remains is to estimate $\Prob[S_1]$.
By definition,
\begin{align}
\Prob[S_1] &= \Prob\left[\,\exists \, i \mid g_i \geq \sigma_0\sqrt{\log n}\right] = 1 - \Prob\left[g_i < \sigma_0\sqrt{\log n}\right]^n\label{eq:Cbound3} \\
&= 1 - \Big(1 - \Prob\big[g_i \geq \sigma_0\sqrt{\log n}\big]\Big)^n.\label{eq:Cbound35}
\end{align}

Using the relationship
$\Phi(x) = \frac{1}{2}\left(1 + \textnormal{erf}\left(\frac{x}{\sqrt{2}}\right)\right)$,
we have
\begin{equation}
\Prob\big[g_i \geq \sigma_0\sqrt{\log n}\big]  \!=\! 1 - \Phi(\sqrt{\log n}) = \!
\frac{1}{2}\!\left(\!\!1 \!-\! \textnormal{erf}\left(\!\!\frac{\sqrt{\log n}}{\sqrt{2}}\right)\!\!\right). 
\label{eq:Cbound45}
\end{equation}
Next, using $\textnormal{erf}(x) \leq \sqrt{1 - \exp\left(-\frac{4}{\pi}x^2\right)}$~\cite{kamath15}
we have
\begin{equation}
\Prob\big[g_i \geq \sigma_0 \sqrt{\log n}\big] \geq \frac{1}{2}\left(1 - \sqrt{1 - \exp\left(-\frac{2}{\pi}\log n\right)}\right).\label{eq:Cbound5}
\end{equation}
A simple calculation shows that
\begin{equation}
\frac{1}{2}\left(1 - \sqrt{1 - \exp\left(-\frac{2}{\pi}\log n\right)}\right) \geq \frac{1}{n}
\label{eq:erfbound}
\end{equation}
if $\frac{\log n}{\log(4n - 4)} \geq \frac{\pi}{2\pi - 2}$, and the first such~$n$ is
$n = 43$. 

Using Equations~\eqref{eq:Cbound5} and~\eqref{eq:erfbound} in Equation~\eqref{eq:Cbound35}, we then find
\begin{equation}
\Prob\big[S_1\big] \geq 1 - \left(1 - \frac{1}{n}\right)^n.
\end{equation}
The right-hand side is decreasing in~$n$, and we take the limit as $n \to \infty$
to derive a bound that holds for all~$n$, giving
$\Prob\big[S_1\big] \geq 1 - e^{-1}$. 
Using this value in Equation~\eqref{eq:Cbound4} gives
\begin{equation}
\mathbb{E}[g_{max}] \geq \big(0.632 - 0.368\sqrt{2/\pi}\big) \sigma_0 \sqrt{\log n},
\end{equation}
and the mean bound follows from a numerical calculation.


The variance bound is proved from the following concentration inequality: defining 
$g = g_{max} - \Exp\left[g_{max}\right]$, 
we have $\Prob(|g|>t) \leq 2 \exp\left(- \frac{t^2}{2\sigma_0^2}\right)$ for $t \geq 0$~\cite[pg. 141]{ChatterjeeSuperconcentration_book:2014}. Define $Z = g^2$, so we have  $\Prob(Z>t) \leq 2 \exp(- \frac{t}{2\sigma_0^2})$ and $\textnormal{var}(X) = \Exp[Z]$.  Since $Z$ is a non-negative r.v.,  
\begin{align*}
\Exp[Z] & = \int_{\R^+}\Prob(Z\geq z)dz  \leq  \int_0^{\infty} 2 \exp(-\frac{z}{2\sigma_o^2})dz \\
& = - 4\sigma_0^2 \left[  \exp\left(-\frac{z}{2\sigma_0^2}\right)  \right]\bigg|^\infty_0 = 4\sigma_0^2. 
\end{align*}
\hfill $\blacksquare$
\end{proof-lemma}

We next apply these results to bound privacy-induced errors in
smart grid calculations. 

\section{Analysis of data analytics accuracy}\label{sec:accuracy-analysis}
This section analyzes the impact of differential privacy first upon customers' billing and
then upon a utility's data-driven analytics that are informed by private data. 
\def\var{\textnormal{var}}
\subsection{Energy Use/Billing}\label{sec:billing}
Recall from Section \ref{sec:problem} that, when using private data, the error in the estimated energy consumption over a period $\K = \{0,1,2,\dots,K\}$ for customer $i$ is
\begin{align}\label{eq:error-billing}
\tilde{E}_\K^i = \Delta t \sum_{k \in \K} w_k^i. 
\end{align}

\paragraph{An arbitrary customer}
We begin our billing analysis by analyzing the impact upon any single customer. 
The mean error in energy consumption
 is $\Exp[\tilde{E}_\K^i ]= 0$, and this holds irrespective of design choices. The variance is
$\var[\tilde{E}_\K^i]  = K[\Delta t]^2\var[w_k^i] = 2K\frac{b^2}{\epsilon^2}[\Delta t]^2$, which
follows from the fact that privacy noise has variance $2c^2 = 2\frac{b^2}{\epsilon^2}$.  
Since $K$ and $\Delta t$ are interrelated, consider a (continuous time) interval $\tau$ over which energy use is to be estimated. Since the corresponding number of  samples is $K = \frac{\tau}{\Delta t}$, we have
\begin{align}\label{eq:var-error-billing}
\var\big[\tilde{E}_\K^i\big] = 2 \tau \Delta t \frac{b^2}{\epsilon^2}.
\end{align}
The trade-off between privacy and accuracy is apparent from this dependency on $b$ and $\epsilon$: a larger $b$ and smaller $\epsilon$ provide stronger privacy (see Section~\ref{sec:privacy}), while Equation~\eqref{eq:var-error-billing} shows that such a choice leads to higher uncertainty in the energy use estimate, with variance growing quadratically
in~$b$ and inverse quadratically with~$\epsilon$. 

A related point is the effect of sampling frequency. Since billing is typically
done monthly, 
we can consider $\tau$ to be a fixed constant and not a design variable. Eq. \eqref{eq:var-error-billing} shows that the uncertainty in the monthly energy use estimate introduced by the privacy noise can be reduced by using a smaller $\Delta t$, i.e., by sampling the demand more frequently\footnote{Although this may appear counterintuitive since the number of random variables being summed increases linearly as $1/\Delta t $, notice that the variance of each summand is quadratic in $\Delta t$, which leads to the reduction in the variance as $\Delta t$ is reduced.}. 

{\bf Cost of privacy:} Let us now consider some numerical values to see how this analysis can drive design of privacy mechanisms. Suppose  $\epsilon=1$ and $\tau = 30 \times 24$ hours (representing one month). For $\Delta t = \frac{1}{4}$ (15-minute sampling), the variance in the monthly energy use is  $360 b^2$ while for $\Delta t = \frac{1}{60}$ (1-minute sampling), the variance reduces to $24b^2$, in (kWh)$^2$. For $b=2$, an arguably small value (as we will see in Section \ref{sec:numerics}), the corresponding standard deviations are  $37.9$ and $9.8$ kWh
for 15-minute and 1-minute sampling, respectively. The average monthly energy use of residential homes in the USA in 2015 was 1883 kWh\footnote{\url{https://www.eia.gov/consumption/residential/data/2015/c\&e/pdf/ce1.1.pdf}}, so a standard deviation of 10 kWh might be tolerable for a consumer, though $38$ kWh might not be.

\paragraph{The worst-affected customer}
Due to the error in the monthly energy use estimate due to privacy noise, a consumer's  monthly bill may also be erroneous. Although these errors are 0 on average, even a single instance of a large error may cause large annoyance, and, in certain  cases, even financial hardship to the customer. The customer that has the largest error will be the one who is most severely affected. The worst error among a set of customers is a crucial value since it may very well dictate consumer acceptance of the technology by driving public debate. We examine this maximum error next. The maximum error is itself a random variable, so its first and second moments are analyzed.


Specializing the discussion following Lemma \ref{lem:minimum-m}  to the sum that makes up $\tilde{E}^i_K$, we see that for $\epsilon = 1$ and $b = 2$, we need $K \geq 22$ for $\tilde{E}^i_K$ to be $0.01$-approximately Gaussian. Even with 15-minute sampling, for a month-long period, $K = 2880 \gg 22$, justifying Gaussian analysis of $\tilde{E}^i_K$.
\begin{theorem}\label{thm:maxerror-billing}
Assuming the Gaussian approximation of $\tilde{E}^i_\K$ holds, the 
	 mean of the maximum of  $\tilde{E}^i_\K$ among $N$ consumers is bounded according to 
	\begin{equation}
      0.478\frac{b}{\epsilon}\Delta t\sqrt{K}\sqrt{\log N}\leq \Exp\left[\max_{1 \leq i \leq N} \tilde{E}^i_\K \right]
	\leq 2\frac{b}{\epsilon}\Delta t\sqrt{K}\sqrt{\log N}.
	\end{equation}
	The variance of the maximum of $\tilde{E}^i_\K$ is upper-bounded by
	\begin{equation}
	\textnormal{var}\left[\max_{1 \leq i \leq N} \tilde{E}^i\right] \leq 8 K\frac{b^2}{\epsilon^2}[\Delta t]^2. 
	\end{equation}
\end{theorem}
\begin{proof-theorem}{\ref{thm:maxerror-billing}}
As shown previously (right before \eqref{eq:var-error-billing}), $\tilde{E}_\K^i$ has mean 0 and variance $2K\frac{b^2}{\epsilon^2}[\Delta t]^2$. The expectation bound then follows
from Lemma~\ref{lem:maxgauss} with $\sigma_0^2 = 2K b^2/\epsilon^2[\Delta t]^2$. The variance bound follows from the second part of Lemma~\ref{lem:maxgauss}. \hfill $\blacksquare$
\end{proof-theorem}

Theorem \ref{thm:maxerror-billing} allows one to estimate the cost of privacy for customers without going through an expensive data collection process. The bounds only depend on statistics of privacy noise, data sampling rate, and number of customers. We will return to this point in Section \ref{sec:numerics}.

\subsection{Aggregate demand}
Recall that $\tilde{d}_k = w_k = \sum_{i\in[N]} w_k^i$ (cf. \eqref{eq:error-aggrdemand}) is the error in the aggregate demand of a collection of $N$ consumers estimated from their private demand data at time~$k$. By appealing to the Berry-Esseen Theorem (Lemma \ref{lem:berry}), we can model $w_k$ as Gaussian, as long as $N$ is large enough. Again specializing the discussion following Lemma \ref{lem:minimum-m}  to the sum that makes up $\tilde{d}_k$, we see that for $\epsilon = 1$ and $b = 2$, we need $N \geq 22$ for $\tilde{d}_k$ to be $0.01$-approximately Gaussian, for any $k$. For any utility, the number of customers is far higher than 22, thus justifying a Gaussian approximation. 

The mean of $\tilde{d}_k$ is 0 and its variance is $\sigma^2_{w_k} = 2N\frac{b^2}{\epsilon^2}$, which follows from the fact that $w_k^i \sim \lap(c)$,  $c = \frac{b}{\epsilon}$, and the variance of a Laplace random variable with scale parameter $c$ is $2c^2$. Therefore, $\Exp[\tilde{d}_k]=0$ and $\var[\tilde{d}_k] = 2N\frac{b^2}{\epsilon^2}$. Thus, though the average error is 0, the uncertainty in the error grows linearly with the number of consumers over which aggregation is performed.


\begin{theorem}\label{thm:maxerror-aggregatedemand}
Assume $N$ is large enough that $\tilde{d}_k$ is approximately Gaussian for all~$k$. Then
	\begin{equation}
	0.478 \sqrt{N}\frac{b}{\epsilon} \sqrt{\log K} \leq 
	\mathbb{E}\left[\max_{1 \leq k \leq K} \tilde{d}_k \right]
	\leq 2\frac{b}{\epsilon}\sqrt{N}\sqrt{\log K}.
	\end{equation}
	The variance of the max, for any $K \in \N$, is bounded via
	\begin{equation}
	\textnormal{var}\left[\max_{1 \leq k \leq K} \tilde{d}_k\right] \leq 8N\frac{b^2}{\epsilon^2}.
	\end{equation}
\end{theorem}
\begin{proof-theorem}{\ref{thm:maxerror-aggregatedemand}}
	Since $\tilde{d}_k = w_k$, which is a sum of~$N$ Laplacian random variables,
	Lemma~\ref{lem:maxgauss} shows that~$N \geq 22$ is sufficient to be
	$0.01$-approximately Gaussian. 
	$w_k$ has mean zero and variance $2Nb^2/\epsilon^2$ 
	and the expectation bound follows from Lemma~\ref{lem:maxgauss} with
	$\sigma_0^2 = 2Nb^2/\epsilon^2$; the variance bound likewise follows
	from Lemma~\ref{lem:maxgauss}. \hfill $\blacksquare$
\end{proof-theorem}

This result shows that on average, the maximum error grows with the square root of the length of time horizon, with variance independent of time. The former is bad news while the latter is good news: peak demand computed over long time horizons may become progressively poorer when it is dominated by the maximum error.  The result also reveals the trade-off between privacy and accuracy: larger $b$ will lead to a larger maximum error, as will smaller values of $\epsilon$.

\section{Numerical results}\label{sec:numerics}
In this section we use high-resolution demand and energy data from a number of residences to illustrate the trade-offs discussed in the previous sections. The data is taken from the Pecan Street Project (from https://dataport.cloud); see~\cite{RhodesExperimental_PecanStreet:2014} for details about the dataset. 

\begin{figure}   
	\includegraphics[width=0.9\columnwidth]{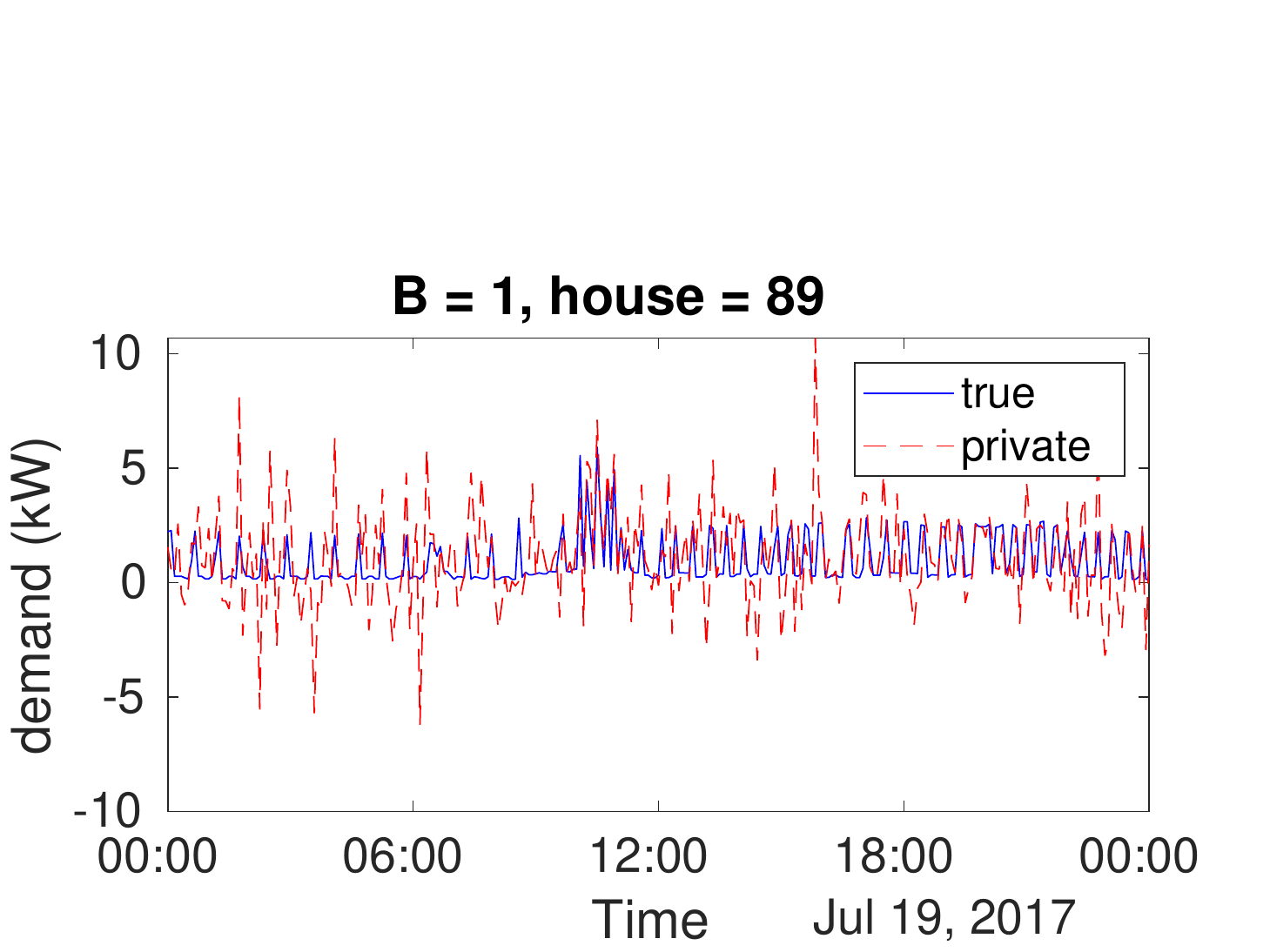}
	\caption{Demand and private demand over 24 hours from an arbitrarily chosen home from the Pecan Street Dataset, with $b=1$. The masking effects of privacy noise can be seen here
	in deviations of the private curve from the true one. 
	}
	\label{fig:demandSingleHomeSingleDay}
\end{figure}

\paragraph{Choosing privacy parameters} For numerical investigations reported here, we choose $\epsilon = \log 2$, which is within typical ranges for differential privacy implementations~\cite{leny14}. The value of $b$ requires more care since it determines what trajectories will be rendered approximately 
indistinguishable with the resulting privacy implementation. Choosing an appropriate $b$ depends on the nature of the data and which events are to be masked~\cite{leny14}. We use $b=1$ for the numerical investigations to follow. This choice creates a significant difference between the true trajectory and the private one. Figure \ref{fig:demandSingleHomeSingleDay} shows  the true demand and its private version created with $b=1$ for an arbitrarily chosen 
home.
We can see from the figure that the privacy noise is larger than even the maximum daily demand in many instants. 



\begin{figure}
\includegraphics[width=0.9\columnwidth]{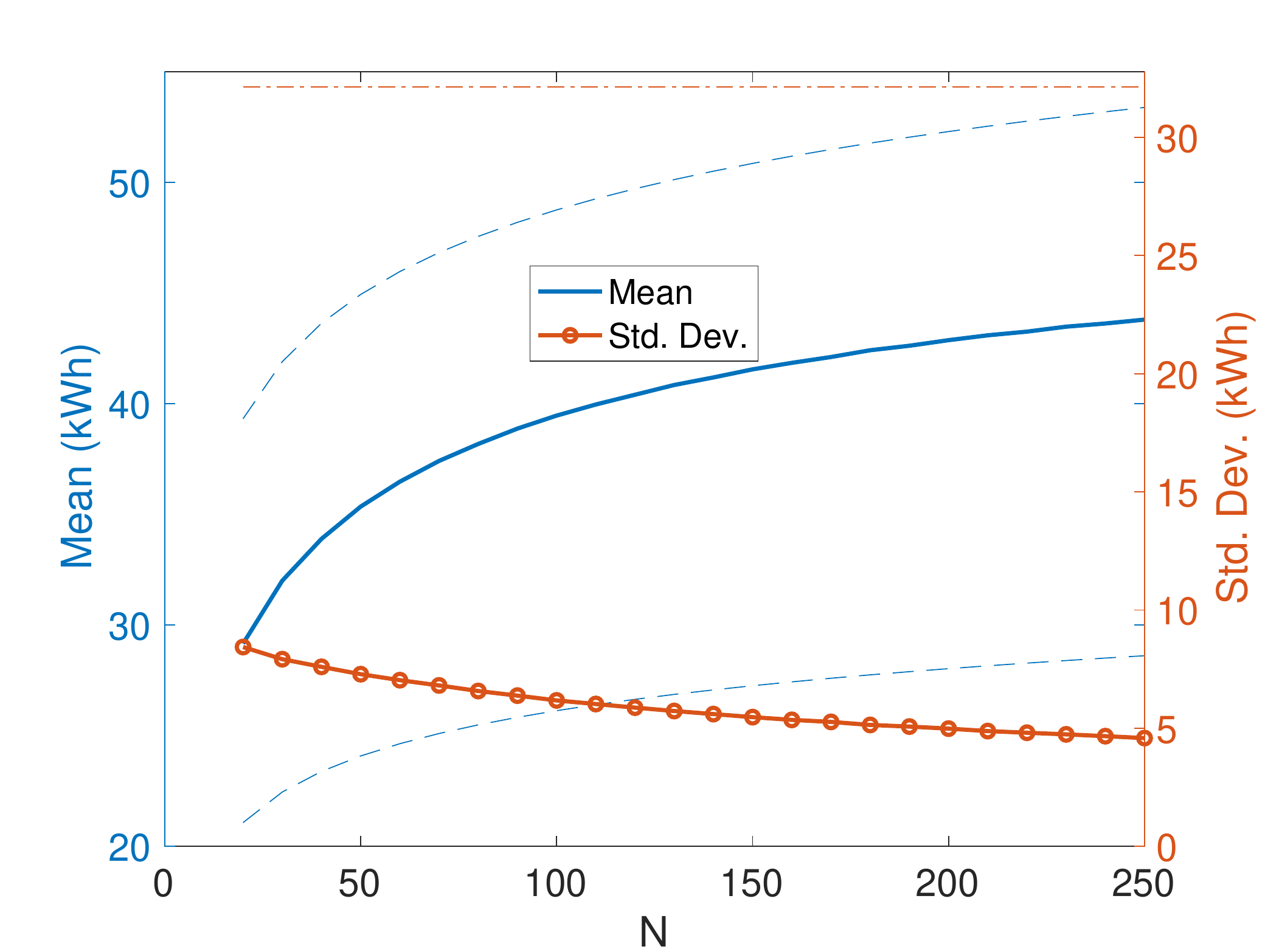}
\caption{Numerically estimated mean and standard deviation of the maximum monthly energy use estimation error as a function of number of customers, computed from 100,000 samples for each $N$. The upper and lower mean bounds from Theorem~\ref{thm:maxerror-billing} are shown in dashed lines, and the upper bound for the standard deviation is shown as a dashdot line.
Theorem~\ref{thm:maxerror-billing} suggests privacy induces only modest error and these
numerical results confirm that this is the case. 
\vspace{-1cm}
}
\label{fig:maxErrorFunctionofNStats}
\end{figure}

\paragraph{Maximum error in monthly energy use} 
Figure \ref{fig:maxErrorFunctionofNStats} shows the numerically estimated mean and variance of the 
maximum 
monthly energy use error 
among $N$ consumers as a function of $N$, for the month of August 2017. Since the demand data is sampled every $5$ minutes, $K = 31 \times 24 \times 12 = 8928$. For each value of $N$ shown, we generate samples of the random variable $Z_{\max}^{(N)} \eqdef \max_{1\leq i\leq N} \tilde{E}^{i}_\K$ via random sampling with replacement, as follows. Among the $300$ total homes we have data for, $N$ homes ($N < 300$) are randomly chosen, the error in their energy use estimate over the month is computed (by using true demand and private demand), and the max value is computed among the $N$ samples. By performing this experiment repeatedly, each time choosing a random subset of size $N$ from the total available $300$, we obtain samples of the random variable $Z_{\max}^{(N)}$. The mean and variance of $Z_{\max}^{(N)}$ are then estimated from 
these samples. 
Figure~\ref{fig:maxErrorFunctionofNStats} also shows the upper and lower bounds on the mean, and the upper bound on the standard deviation, predicted by Theorem~\ref{thm:maxerror-billing}. We see the bounds on the mean are tighter than the bound on the variance. 
 
{\bf Cost of privacy (again):}  As in the single/arbitrary customer case, the maximum error in monthly energy use estimate due to privacy can be used to compute a worst-case cost of privacy. In this instance, for a utility with a customer base of 250 households, the average cost of privacy to the worst-hit customer is $43.8$kWh/month, which translates to \$$5.6$/month at the current average rate of $12.9 \cent$/kWh\footnote{From \url{https://www.eia.gov/energyexplained/}, second tab.}.  If one uses ``mean plus 3-sigma'' to estimate the cost, then it turns out to be $43.8+15 = 58.8$ kWh/month, or \$$7.6$/month. For a larger customer base, these numbers will increase, but since the growth of the mean is logarithmic and the variance is constant, the increase will be small. 
  These bounds allow us to compute the same worst-case cost of privacy but without having to go through an expensive data collection process on customer demand. Repeating the calculation done in the previous paragraph, but using the theoretical upper bounds obtained in Theorem \ref{thm:maxerror-billing}, we find an upper bound on the cost of privacy for the worst-hit customer to be \$$19$/month. 

\paragraph{Maximum error in aggregate demand} 

\begin{figure}
	\includegraphics[width=0.9\columnwidth]{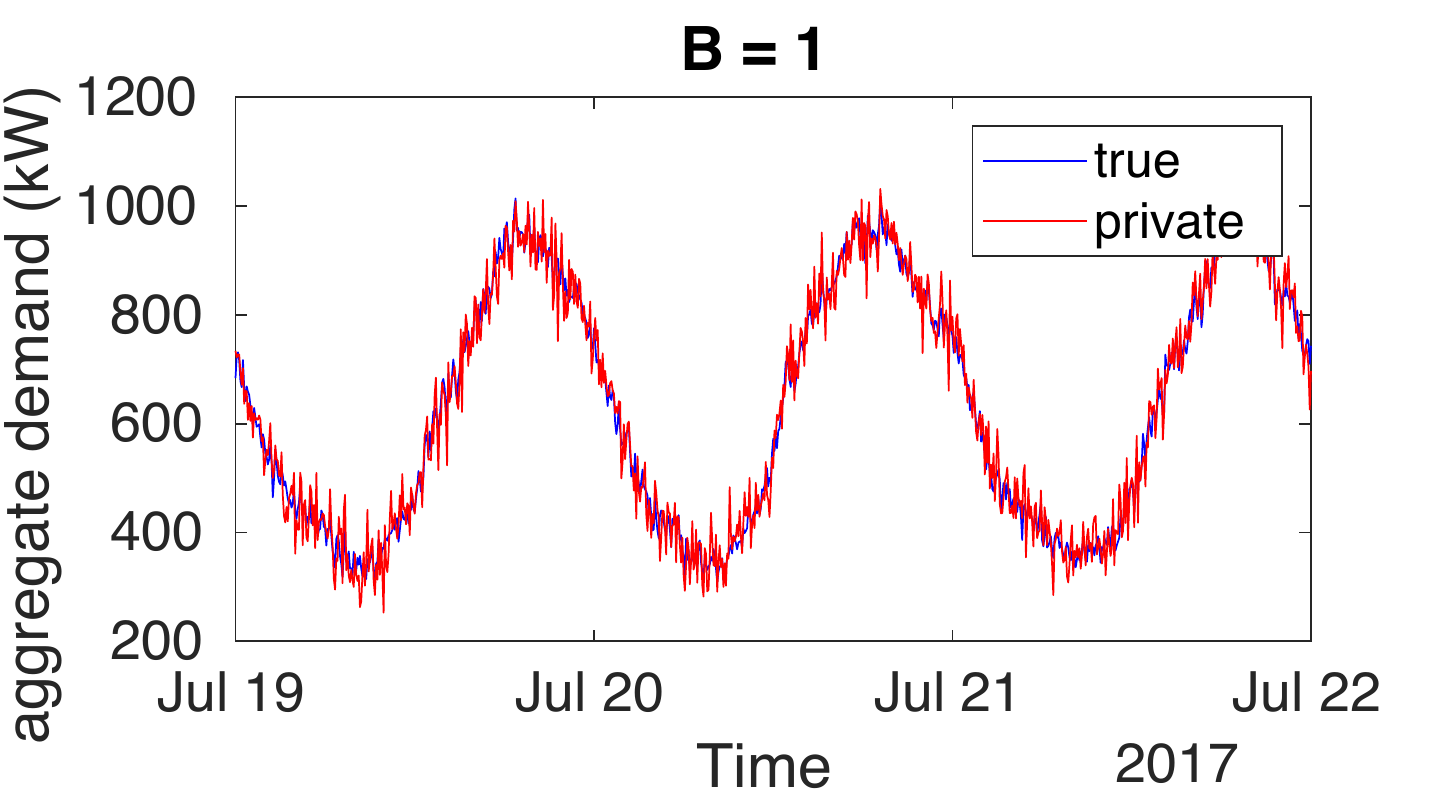}
	\caption{Aggregate demand and private aggregate demand from 300 homes in the Pecan Street Project dataset over 4 days. Here we see only modest error in the private aggregate
	demand signal, indicating that accurate data can be gathered even under a mandate
	of privacy. 
	}
	\label{fig:aggregateDemand}
\end{figure}

\begin{figure}
\includegraphics[width=0.9\columnwidth]{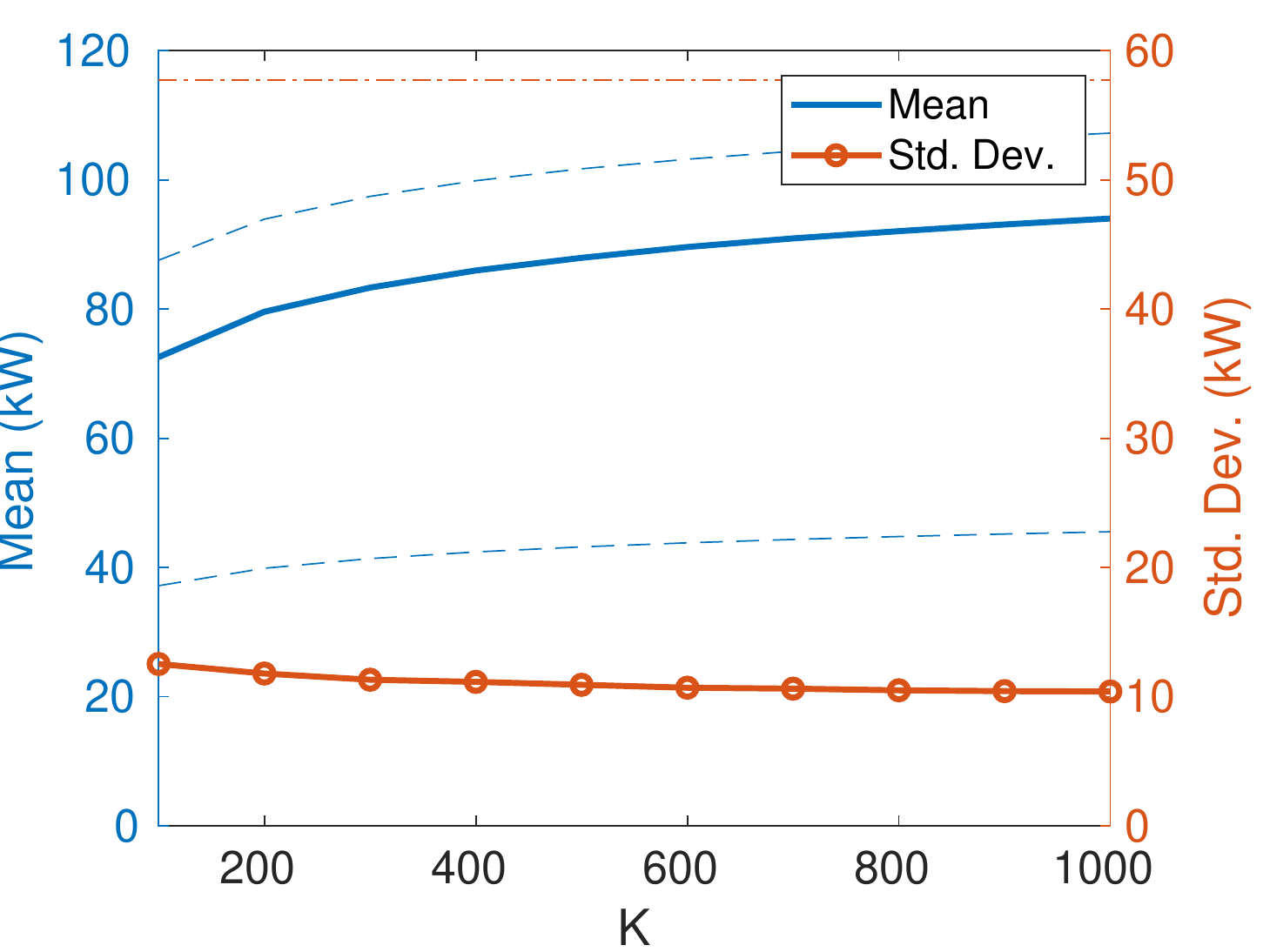}
\caption{
Numerically estimated mean and standard deviation of the maximum monthly energy use estimation error as a function of times, computed from data for 300 homes. The upper and lower mean bounds from Theorem~\ref{thm:maxerror-aggregatedemand} are shown in dashed lines, and the upper bound for the standard deviation is shown as a dashdot line.
Theorem~\ref{thm:maxerror-aggregatedemand} suggests privacy induces only modest error and these
numerical results confirm that this is the case.  
}
\label{fig:maxerrorK}
\end{figure}

Figure \ref{fig:aggregateDemand} shows the trajectory of the aggregate demand from all 300 homes in the dataset, as well as the private aggregate demand (computed by using the private demand data), for a period of 4 days. Figure~\ref{fig:maxerrorK} shows mean and standard deviation bounds on the maximum
error in the private data over time. Both plots indicate that privacy incurs only
modest errors with respect to power usage data. 

\section{Conclusion}\label{sec:conclusion}
We applied trajectory-level differential privacy to demand data from customer smart meters. Because data is made private at the customer's home, any analysis with that data has uncertainty. This includes computing monthly energy use, which is essential for billing. We analyzed the average and worst case errors, and showed that the trade-off between privacy guarantees and analysis accuracy can be translated to a financial cost of privacy to the consumer. 

In this preliminary work, we limited our analysis to energy use for a consumer and aggregate demand among all consumers. The peak demand of the aggregate is a valuable quantity for grid planning, and how privacy noise affects the estimate of peak demand remains to be investigated.


\bibliographystyle{IEEEtran}
\bibliography{sources}
\end{document}